\input amstex.tex

\documentstyle{amsppt}
\magnification=1200

\loadbold

\def\Co#1{{\Cal O}_{#1}}
\def\Fi#1{\Phi_{|#1|}}
\def\fei#1{\phi_{#1}}
\def\rw{\rightarrow}
\def\lrw{\longrightarrow}
\def\KX{K_X}

\def\Bbbp1{{\Bbb P}^1}
\def\simlin{\sim_{\text{\rm lin}}}
\def\simnum{\sim_{\text{\rm num}}}
\def\Ce#1{{\Cal E}_{#1}}

\def\Div{\text{\rm Div}}

\TagsOnRight
\pagewidth{13cm}
\pageheight{19.5cm}
\vcorrection{-0.3in}
\topmatter
\title Kawamata-Viehweg Vanishing and Quint-canonical\\ Map 
of a Complex Threefold
\endtitle
\author \smc Meng Chen$^*$\endauthor
\thanks
The author is partially supported by  the National Natural Science Foundation
of China
\endthanks
\endtopmatter
\document
\baselineskip 14pt
\leftheadtext{Chen}
\rightheadtext{quint-canonical map}
\centerline{\it Department of Applied Mathematics}
\centerline{\it Tongji University, Shanghai, 200092, China}
\vskip0.2cm

\centerline{\it e-mail: chenmb\@online.sh.cn}
\vskip0.3cm

\head Introduction \endhead
Given a complex nonsingular minimal threefold $X$ of general type, Benveniste
(\cite{2}) proved that m-canonical map $\phi_m$ is a birational map onto
its image when $m\ge 8$, Matsuki (\cite{15}) showed the same statement for 
7-canonical map. In \cite{5}, we proved the birationality of 6-canonical map. 
In \cite{14}, Lee  proved,
 independently, that m-canonical map is a birational morphism for $m\ge 6$.
Furthermore, the 5-canonical map is birational when $K_X^3>2$ according to
Ein-Lazarsfeld-Lee. 
The aim of this note is to prove the following two theorems
by a different method:

\proclaim{\smc Theorem 1} Let $X$ be a complex nonsingular projective threefold
with nef and big canonical divisor $K_X$. Then

(1) $\fei{5}$ is a birational map onto its image when $p_g(X)\ge 3$;

(2) if $p_g(X)=2$ and $\fei{5}$ is not a birational map, then $\fei{5}$ is  generically
finite of degree 2 and $q(X)=h^2(\Co{X})=0$ and $|K_X|$ is composed of a rational pencil
of surfaces of general type with $(K^2,p_g)=(1,2)$.
\endproclaim

\proclaim{\smc Theorem 2} Let $X$ be a complex nonsingular projective threefold
with nef and big canonical divisor $K_X$. Suppose $p_g(X)\le 1$ and
$|2K_X|$ be composed of a pencil of surfaces, i.e., $\text{dim}\phi_2(X)=1$, 
then $\phi_5$ is a birational map onto its image.
\endproclaim

We would like to put a conjecture here:
\proclaim{\smc Conjecture} There exists the exception and 
the only possible exception to the birationality
of 5-canonical map of a complex nonsingular minimal threefold $X$ is one
with 
$$(K_X^3, p_g(X), q(X), h^2(\Co{X}))=(2, 2, 0, 0).$$
\endproclaim

\head 1. A lemma on a surface with $K^2=1$ and $p_g=2$ \endhead

Let $S$, with minimal model $S_0$,
 be a nonsingular algebraic surface of general type with
$K_{S_0}^2=1$ and $p_g(S)=2$. It is well-known that $\fei{5}$ is birational
and $\fei{4}$ is generically finite of degree 2. In order to make  preparation
for the proof of our main theorems. We would like to formulate a remark
to this kind of surfaces.

Kawamata-Viehweg vanishing theorem will be used throughout this paper in the following
form:
\proclaim{\smc Vanishing Theorem} Let $X$ be a nonsingular complete variety,
$D\in\Div(X)\otimes {\Bbb Q}$. Assume the following two conditions:

(1) $D$ is nef and big;

(2) the fractional part of $D$ has the support with only normal crossings.

Then $H^i(X, \Co{X}(\lceil D\rceil+K_X))=0$ for $i>0$, where $\lceil D\rceil$ 
is the minimum integral divisor with $\lceil D\rceil-D\ge 0$.
\endproclaim

\remark{\smc Remark 1.1} In the case of surfaces, Sakai proved that the
Kawamata-Viehweg vanishing holds without the assumption of normal crossings.
\endremark

\proclaim{\smc Lemma 1.1} Let $S$, with minimal model $S_0$,
be a nonsingular projective algebraic
surface of general type with $K_{S_0}^2=1$ and $p_g(S)=2$. If $\pi: S\lrw S_0$
is the contraction map, then
$$\phi_{4.5}:=\Fi{K_S+3\pi^*(K_{S_0})+\lceil\frac{\pi^*(K_{S_0})}{2}\rceil}$$
is a birational map onto its image.
\endproclaim
\demo{\smc Proof}
If $\pi^*(K_{S_0})$ is an irreducible effective divisor, the lemma is obviously
true. Otherwise, we have an effective irreducible divisor $D_0$ and an 
effective divisor $E_0$ such that $D_0+E_0\in |\pi^*(K_{S_0})|$ and 
$D_0\cdot \pi^*(K_{S_0})=1$.

We know that $|K_{S_0}|$ has exactly one base point and has no fixed part.
(one may consult (8.1) at page 225 of \cite{1}). A general member $C\in
|K_{S_0}|$ is a nonsingular curve of genus 2. Let $P$ be the base point of
$|\pi^*(K_{S_0})|$.
It is obvious that 
$\Fi{K_S+3\pi^*(K_{S_0})+\lceil\frac{D_0+E_0}{2}\rceil}$ can separate 
two general members of $|\pi^*(K_{S_0})|$. We may suppose $S$ be like one of 
the following three cases without losing of generality:

(1) the exceptional divisors of $\pi$ do not lie over $P$;

(2) $S$ is just obtained by blowing up the base point $P$ from a surface 
like case (1);

(3) $S$ is obtained by several blow ups from a surface like case (2).

Case (1). In this case, there is no changes around $P$. So we again denote
$\pi^{-1}(P)$ by $P$ with no confusion. Let $\tilde C$ 
be the strict transforms of $C$. Denote $\overline{D_0}=\pi_*{D_0}$
and $\overline{E_0}=\pi_*(E_0)$.

Let $K_S=\pi^*(K_{S_0})+\sum E_j$. Note that $\pi |_{\tilde{C}}:\tilde{C}
\lrw C$ is an isomorphism. 
Because $3\pi^*(K_{S_0})+\frac{D_0+E_0}{2}-\tilde{C}\simnum 
\frac{5}{2}\pi^*(K_{S_0})$ 
is nef and big, therefore, by Vanishing Theorem, we have
$$H^1(S, K_S+3\pi^*(K_{S_0})+\lceil\frac{D_0+E_0}{2}\rceil-\tilde{C})=0.$$
Note that, in this case, $K_S|_{\tilde{C}}=\pi^*(K_{S_0})|_{\tilde{C}}$
and $\tilde{C}\in|\pi^*(K_{S_0})|$.
We see that 
$$\Fi{K_S+3\pi^*(K_{S_0})+\lceil\frac{D_0+E_0}{2}\rceil}|_{\tilde{C}}
=\Fi{2K_{\tilde{C}}+q},$$
where $q:=D_0|_{\tilde{C}}$ is a point on $\tilde{C}$. Because
$\deg(2K_{\tilde{C}}+q)=5$ and then
$2K_{\tilde{C}}+q$ is very ample, 
$$\Fi{K_S+3\pi^*(K_{S_0})+
\lceil\frac{D_0+E_0}{2}\rceil}$$ 
is a birational map onto its image.

Case (2). In this case, let $S_1$ be a surface as case (1) and $\pi_1:S_1\lrw S_0$
be the contraction map onto $S_0$. Let $\pi_2:S\lrw S_1$ be the blowing up
at $P$, i.e., the base point of $\pi_1^*(K_{S_0})$. Let $C\in
|\pi_1^*(K_{S_0})|$ be a general member and $\tilde{C}$ the strict transform
of $C$. Let $D_1:={\pi_1}_*D_0$ and $E$ be the $(-1)$-curve over $P$.
Denote $\pi:=\pi_1\circ\pi_2$.

We have
$$\align
K_S&=\pi_2^*(K_{S_1})+E\\
   &=\pi_2^*(\pi_1^*(K_{S_0})+\sum E_k)+E\\
   &=\pi^*(K_{S_0})+\pi_2^*(\sum E_k)+E.
\endalign$$
We also have that $\pi^*(K_{S_0})\simlin \tilde{C}+E$.

Now we consider the system 
$$|K_S+2\pi^*(K_{S_0})+\tilde{C}+\lceil\frac{D_0+E_0}{2}\rceil|.$$
Because 
$$K_S+2\pi^*(K_{S_0})+\tilde{C}+\lceil\frac{D_0+E_0}{2}\rceil\le
K_S+3\pi^*(K_{S_0})+\lceil\frac{D_0+E_0}{2}\rceil,$$
we only have to verify the birationality of
$$\Fi{K_S+2\pi^*(K_{S_0})+\tilde{C}+\lceil\frac{D_0+E_0}{2}\rceil}.$$
Because 
$2\pi^*(K_{S_0})+\frac{D_0+E_0}{2}$ is nef and big,
we have 
$$H^1(S,K_S+2\pi^*(K_{S_0})+\lceil\frac{D_0+E_0}{2}\rceil)=0$$
by Vanishing Theorem. Note that $E_0=E+E'$, $E'\ge 0$ and $2E\not\le E_0$.
Therefore 
$$\Fi{K_S+2\pi^*(K_{S_0})+\tilde{C}+\lceil\frac{D_0+E_0}{2}\rceil}|_{\tilde{C}}
=\Fi{2K_{\tilde C}+q},$$
where $q=E|_{\tilde{c}}$. $\Fi{2K_S+q}$ is an embedding, because $\deg
(2K_S+q)=5$. Thus 
$$\Fi{K_S+3\pi^*(K_{S_0})+\lceil\frac{D_0+E_0}{2}\rceil}$$
is birational.

Case (3). one can easily go through the proof by a similar argument as that
of case (2).
\qed\enddemo

\head 2. Proof of theorem 1\endhead
\subhead Basic formula\endsubhead
Let $X$ be a nonsingular projective threefold. For a divisor $D\in \Div(X)$, 
we have
$$\chi(\Co{X}(D))=D^3/6-K_X\cdot D^2/4+D\cdot(K_X^2+c_2)/12+\chi(\Co{X})$$
by Riemann-Roch theorem. A calculation shows that
$$\chi(\Co{X}(D))+\chi(\Co{X}(-D))=-K_X\cdot
D^2/2+2\chi(\Co{X})\in{\Bbb Z},$$
therefore $K_X\cdot D^2$ is an even integer, especially $K_X^3$
is even.
If $\KX$ is nef and big, then we obtain by Kawamata-Viehweg's vanishing
theorem that
$$p(n):=h^0(X,\Co{X}(n\KX))=(2n-1)[n(n-1)\KX^3/12-\chi(\Co{X})],\tag 2.1$$
for $n\ge 2$. 
 Miyaoka (\cite{16})  showed that $3c_2(X)-c_1(X)^2$ is
pseudo-effective, therefore we get $\KX^3\le -72\chi(\Co{X})$ by
the Riemann-Roch equality, $\chi(\Co{X})=-c_2\cdot\KX/24.$
In particular, $\chi(\Co{X})<0$.

Let $f:X\rw C$ be a fibration onto a nonsingular
curve $C$. From the spectral sequence:
$$E_2^{p,q}:=H^p(C,R^qf_*\omega_X)\Longrightarrow
E^n:=H^n(X,\omega_X),$$
a direct calculation shows that
$$h^2({\Cal O}_X)=h^1(C,f_*\omega_X)+h^0(C,R^1f_*\omega_X), \tag 2.2$$
$$q(X):=h^1({\Cal O}_X)=b+h^1(C,R^1f_*\omega_X). \tag 2.3 $$
Therefore we obtain
$$\chi({\Cal O}_X)=\chi({\Cal O}_F)\chi({\Cal O}_C)+\Delta_2-\Delta_1, \tag 2.4$$
where we set $\Delta_1:=\deg f_*\omega_{X/C}$ and
$\Delta_2:=\deg R^1f_*\omega_{X/C}$.  Theorem 1 of \cite{11} tells that
$\Delta_1\ge 0$. Lemma 2.5 of \cite{17} says that $\Delta_2\ge 0$.

\proclaim{\smc Lemma 2.1} Let $S$ be a nonsingular algebraic surface, $L$ a nef
and big divisor on $S$. Then

(1) $\Fi{K_S+mL}$ is a birational map onto its image for $m\ge 4$;

(2) $\Fi{K_S+3L}$ is a birational map onto its image when $ L^2\ge 2$.
\endproclaim
\demo{\smc Proof} This is a direct result of Corollary 2 of \cite{18}.
\qed\enddemo

\proclaim{\smc Lemma 2.2} (See Lemma 2 of \cite{19})
Let $X$ be a nonsingular projective variety, $D$ 
a divisor with $|D|\ne \emptyset$.
If the complete linear system $|M|$ is base point free and $\text{\rm dim}
\Fi{M}(X)\ge 2$ and $\Fi{M+D}$ is not a birational map onto its image,
then $\Fi{M+D}|_S$ is also not birational for a general member $S\in |M|$.
\endproclaim

\definition{\smc Definition 2.1} Let $X$ be a nonsingular projective threefold.
If $\dim\fei{1}(X)\ge 2$ and set $K_X\simlin M_1+Z_1$, where $M_1$ is
the moving part and $Z_1$ the fixed one. We define $\delta_1(X):=K_X^2\cdot M_1$.
\enddefinition

\proclaim{\smc Proposition 2.1} Let $X$ be a nonsingular projective threefold
with nef and big canonical divisor $K_X$. Suppose  
$\dim\fei{1}(X)\ge 2$, then $\delta_1(X)\ge 2$.
\endproclaim
\demo{\smc Proof}
Let $f_1:X'\lrw X$ be a succession of blowing-ups with nonsingular centers according
to Hironaka such that $g_1:=\fei{1}\circ f_1$ is a morphism. Let
$g_1:X' @> h_2>>W_1'@> s_1>>{W_1\subset{\Bbb P}^{p_g(X)-1}}$
be  the Stein factorization of $g_1$. Let $H_1$ be a
hyperplane section of $W_1=\overline{\fei{1}(X)}$ in ${\Bbb P}^{p_g(X)-1}$ 
and $S_1$ be a general member of $|g_1^*(H_1)|$.
Since $\dim W_1\ge 2$, $S_1$ is a nonsingular irreducible
projective surface by Bertini Theorem.   
Set $f_1^*(M_1)\simlin S_1+E_1'$,
$K_{X'}\simlin f_1^*(K_X)+E_1$, where
$E_1$ is the ramification divisor for $f_1$, $E_1'$ is the
exceptional divisor for $f_1$. We have the following commutative diagram:
$$\CD
X'   @>h_1>>    W_1'\\
@|            @VVs_1V\\
X'   @>>g_1>   W_1\\
@Vf_1VV            @.\\
X    @.
\endCD$$

We have $\delta_1(X)=K_X^2\cdot M_1=f_1^*(K_X)^2\cdot S_1$.
Multiplying $K_X\simlin M_1+Z_1$ by
$K_X\cdot M_1$, we have
$$K_X^2\cdot M_1=K_X\cdot M_1^2+K_X\cdot M_1\cdot Z_1.$$
Since $|S_1|$ is not composed of a pencil, $f_1^*(K_X)$ is nef
and big and since $S_1$ is nef, we have $f_1^*(K_X)\cdot S_1^2
\ge 1$. So that
$$\align
K_X\cdot M_1^2&=f_1^*(K_X)\cdot f_1^*(M_1)^2=f_1^*(K_X)\cdot
f_1^*(M_1)\cdot S_1\\
&=f_1^*(K_X)\cdot S_1^2+f_1^*(K_X)\cdot S_1\cdot E_1'\ge 1.
\endalign $$
Whereas, $K_X\cdot M_1^2$ is even  and
$K_X\cdot M_1\cdot Z_1\ge 0$ because $M_1\cdot Z_1\ge 0$ as a
1-cycle. Thus we have $K_X^2\cdot M_1\ge 2$. 
\qed\enddemo

\proclaim{\smc Theorem 2.1} Let $X$ be a nonsingular projective complex 
threefold with nef and big canonical divisor $K_X$. Suppose $p_g(X)\ge 3$
and $|K_X|$ be not composed of a pencil of surfaces, i.e., $\dim\fei{1}(X)\ge 2$,
then $\fei{5}$ is a birational map onto its image.
\endproclaim
\demo{\smc Proof} 
We use the same diagram as in the proof of Proposition 2.1 and keep the same 
notations there. Assume $\fei{5}$ be not birational, because
$$5K_{X'}\simlin \{K_{X'}+3f_1^*(K_X)+S_1\}+4E_1+f_1^*(Z_1)+E_1',$$
$\Fi{K_{X'}+3f_1^*(K_X)+S_1}$ is also not birational. Therefore 
$\Fi{K_{X'}+3f_1^*(K_X)+S_1}|_{S_1}$ is not birational by Lemma 2.2.

On the other hand, we have $H^1(X', K_{X'}+3f_1^*(K_X))=0$ according to 
Vanishing Theorem. Thus 
$$\Fi{K_{X'}+3f_1^*(K_X)+S_1}|_{S_1}=\Fi{K_{S_1}+3L_1},$$
where we set $L_1:=f_1^*(K_X)|_{S_1}$, which is nef and big and $L_1^2
=\delta_1(X)\ge 2$. Therefore the latter is birational onto its image
by Lemma 2.1. Which is a contradiction.
\qed\enddemo

In the next, we always suppose that $|K_X|$ be composed of a pencil of surfaces.
We again use the same diagram as in the proof of Proposition 2.1. Note 
that $W_1'$ is a nonsingular curve. We usually call $h_1$ a derived fibration
of $\fei{1}$. Let $F$ be a general fiber of $h_1$. Then $F$ must be a nonsingular 
projective surface of general type by Bertini Theorem. Denote $b:=g(W_1')$, the 
geometric genus of curve $W_1'$. 

We can set $g_1^*(H_1)\simnum aF,$ where $a\ge p_g(X)-1$. Let $\overline{F}:={f_1}_*(F)$,
then $M_1\simnum a\overline{F}$. We will formulate our proof through two steps:
(1) $K_X\cdot \overline{F}^2>0$ and (2) $K_X\cdot \overline{F}^2=0$.

\proclaim{\smc Theorem 2.2} Let $X$ be a nonsingular projective threefold
with nef and big canonical divisor $K_X$. Suppose $p_g(X)\ge 2$ and 
$|K_X|$ be composed of a pencil, keeping the above notations,  
if $K_X\cdot \overline{F}^2>0$,
then $\fei{5}$ is a birational map onto its image.
\endproclaim
\demo{\smc Proof}
We have 
$$5K_{X'}\simlin \{K_{X'}+3f_1^*(K_X)+aF\}+4E_1+f_1^*(Z_1)+E_1'.$$
Consider the system $|K_{X'}+3f_1^*(K_X)+aF|$, we have
$H^1(X', K_{X'}+3f_1^*(K_X))=0$. Generically, we can take $g_1^*(H_1)$
be a disjoint union of fibers $F_i$ ($1\le i\le a$). Therefore we have
the following exact sequence:
$$\align
0\lrw\Co{X'}(K_{X'}+3f_1^*(K_X))&\lrw\Co{X'}(K_{X'}+3f_1^*(K_X)+
  g_1^*(H_1))\\
  &\lrw\oplus_{i=1}^a\Co{F_i}(K_{F_i}+3L_i)\lrw 0,
\endalign$$
where $L_i=f_1^*(K_X)|_{F_i}$, which is nef and big and
$$L_i^2=K_X^2\cdot \overline{F}\ge K_X\cdot\overline{F}^2\ge 2.\ 
(K_X\cdot\overline{F}^2 \ \text{\rm is even})    $$
From the above exact sequence, we see that
$$\Fi{K_{X'}+3f_1^*(K_X)+g_1^*(H_1)}|_{F_i}=\Fi{K_{F_i}+3L_i}$$
is a birational map onto its image by Lemma 2.1. Thus $\fei{5}$ is 
birational.
\qed\enddemo

\proclaim{\smc Lemma 2.3} Let $X$ be a nonsingular projective threefold
with nef and big canonical divisor $K_X$. Keeping the above notations,
if $K_X\cdot\overline{F}^2=0$, then
$$\Co{F}(f_1^*(K_X)|_F)\cong\Co{F}(\pi^*(K_{F_0})).$$
\endproclaim
\demo{\smc Proof} This can be obtained by a similar argument to that
for {\it Case} $\beta)$ of $(ii)$, Theorem 7 of \cite{15}.
\qed\enddemo

\proclaim{\smc Theorem 2.3} Under the same assumption as in Lemma 2.3,
If the minimal model $F_0$ of $F$ is not a surface with $K_{F_0}^2=1$
and $p_g(F_0)=2$, then $\fei{5}$ is a birational map onto its image.
\endproclaim
\demo{\smc Proof} The proof is almost the same as that of Theorem 2.2.
The only difference occurs on $L_i=f_1^*(K_X)|_{F_i}$. From Lemma 2.3,
we see that $L_i\simlin \pi^*(K_{F_0})$ and therefore
$\Fi{K_{F_i}+3L_i}=\Fi{4K_{F_i}}$ is birational under the assumption
of this theorem.
\qed\enddemo

\proclaim{\smc Theorem 2.4} Under the same assumption as in Lemma 2.3,
if the minimal model $F_0$ of $F$ is just the surface with $K_{F_0}^2=1$
and $p_g(F_0)=2$, then $\fei{5}$ is also a birational map in one of the 
following two cases:

(1) $p_g(X)\ge 3$;

(2) $ p_g(X)=2$ and $b:=g(W_1')\ne 0$.
\endproclaim
\demo{\smc Proof} In the two cases of this theorem, we can see that $a\ge 2$. 
Fix an effective divisor $K_0\in|K_X|$. Actually, we can modify $f_1$
such that 
$$f_1^*(K_0)=\sum_{i=1}^aF_i+E_1'+f_1^*(Z_1)$$
has support with only normal crossings. 
Thus, from now on, we always suppose $f_1$ has this property.
For a general fiber $F$ of $h_1$, 
We have $g_1^*(H_1)\simnum 2F+\sum_{i=1}^{a-2}F_i$.
For ${\Bbb Q}$-divisor 
$$\overline{G}:=4f_1^*(K_X)-F-\frac{1}{2}(F_1+\cdots+F_{a-2}
+E_1'+f_1^*(Z_1)),$$
 it is nef and big. Denote 
$$G:=[\frac{F_1+\cdots+F_{a-2}+E_1'+f_1^*(Z_1)}{2}],$$
then $H^1(X', K_{X'}+4f_1^*(K_X)-F-G)=0$ by Vanishing Theorem.
 Considering the system $|K_{X'}+4f_1^*(K_X)-G|$, it is obvious that
$$K_{X'}+4f_1^*(K_X)-G\le 5K_{X'}.$$ 
In order to proof the birationality of $\fei{5}$, we only have to
verify for 
$$\Fi{K_{X'}+4f_1^*(K_X)-G}.$$ 
From the following exact sequence
$$\align
0&\lrw\Co{X'}(K_{X'}+4f_1^*(K_X)-G-F)\lrw\Co{X'}(K_{X'}+4f_1^*(K_X)-G)\\
  &\lrw\Co{F}(K_F+3f_1^*(K_X)|_F+\lceil\frac{ E_1'|_F+f_1^*(Z_1)|_F}{2}\rceil)
\lrw 0,
\endalign$$
we see that 
$$\Fi{K_{X'}+4f_1^*(K_X)-G}|_F=\Fi{K_F+3f_1^*(K_X)|_F+\lceil\frac{E_1'|_F+
f_1^*(Z_1)|_F}{2}\rceil}.$$
Note that $f_1^*(K_X)|_F\simlin E_1'|_F+f_1^*(Z_1)|_F$. 
From Lemma 2.3, we have 
$f_1^*(K_X)|_F\simlin\pi^*(K_{F_0})$, where $\pi:F\lrw F_0$ is the contraction
to the minimal model. 
 Thus we complete the proof 
by Lemma 1.1.
\qed\enddemo

Finally, if $p_g(X)=2$ and $|K_X|$ is composed of a pencil of surfaces,
 the above method is not effective.
But from the proof of Theorem 2.3, we can see that $\fei{5}$ is at least
a generically finite map of degree $2$. 
By formula (2.2) and (2.3), we can easily get $q(X)=h^2(\Co{X})=0$.

Combining the arguments of this section, we obtain Theorem 1.

\head 3. On a bicanonical pencil of surfaces of general type\endhead
In order to study the case when $p_g(X)\le 1$, it is natural to study
$\fei{2}$. This section is a preparation for the proof of Theorem 2.

Let $X$ be a nonsingular minimal projective threefold. If $|2\KX|$ 
is composed of a pencil of surfaces, i.e., the image of $X$ through $\Fi{2K_X}$
is of dimension 1, we can find a birational modification
$f_2:X'\rw X$ such that $g_2=\Fi{2\KX}\circ f_2$ is a morphism. Let 
 $W_2=\overline{\phi_2(X)}\subset{\Bbb P}^{p(2)-1}$, 
and $g_2=s_2\circ h_2$ is a Stein-factorization of $g_2$. 
We have the following commutative diagram:
$$\CD
X'   @>h_2>>    C\\
@|            @VVs_2V\\
X'   @>>g_2>   W_2\\
@Vf_2VV            @.\\
X    @.
\endCD$$
where $h_2:X'\rw C$ is called a derived fibration of $\phi_2$.
Let $F$ be a general fiber of $h_2$, 
then $F$ must be a nonsingualar projective surface by Bertini Theorem.
Denote  $b:=g(C)$, the genus of $C$.

\proclaim{\smc Lemma  3.1} (Claim 9.1 of \cite{15})
Let $X$ be a nonsingular minimal projective threefold of general type,
if $|2\KX|$ is composed of a pencil of surfaces, then
$$\Co{F}(f_2^*(\KX)|_F)\cong\Co{F}(\pi^*(K_{F_0})),$$
where  $\pi:F\rw F_0$ is the birational contraction onto the minimal model.
\endproclaim
\proclaim{\smc Lemma 3.2} Under the same assumption as in Lemma 3.1,
we have $K_{F_0}^2\le 3$ and $F$ is of one of the following two cases:
\roster
\item $q(F)=0$, $p_g(F)\le 3$;
\item $p_g(F)=q(F)=1$.
\endroster
\endproclaim
\demo{\smc Proof} 
Let $f_2^*(2\KX)\simlin g_2^*(H_2)+Z_2'$, where  $Z_2'$ is the fixed part
and $H_2$ is a general hyperplane section of $W_2$. Obviously we have
$g_2^*(H_2)\simnum a_2F$, $a_2\ge p(2)-1$. From Lemma 3.1, we have
$$K_{F_0}^2=(f_2^*(\KX)|_F)^2=f_2^*(\KX)^2\cdot F.$$
Let $2\KX\simlin M_2+Z_2$, where $M_2$ is the moving part and $Z_2$ is
 the fixed part. We also have
$M_2={f_2}_*(g_2^*(H_2))$. Denote $\overline{F}={f_2}_*F$, then  
$M_2\simnum a_2\overline{F}$.
By projection formula, one has
$$K_X^2\cdot \overline{F}=f_2^*(\KX)^2\cdot F=K_{F_0}^2.$$
Because $K_X$ is  nef, we have
$2K_X^3\ge a_2K_X^2\cdot\overline{F}.$ Therefore
$$K_X^2\cdot\overline{F}\le\frac{2}{a_2}K_X^3\le\frac{4K_X^3}{K_X^3-6\chi(\Co{X})-2}
\le\frac{4K_X^3}{K_X^3+4}<4,$$
and then $K_{F_0}^2\le 3$.
Because $2p_g(F_0)-4\le K_{F_0}^2$,  $p_g(F_0)\le 3$.
If $q(F)>0$, then Bombieri's theorem(\cite{3}) tells that  $K_{F_0}^2\ge 2\chi(\Co{F_0})\ge
2$, therefore $\chi(\Co{F_0})=1$, i.e., $p_g(F_0)=q(F_0)$. 
By  Debarre's result(\cite{7}), we have   $K_{F_0}^2\ge 2p_g(F_0)$, 
therefore $p_g(F_0)=1$. \qed\enddemo

\proclaim{\smc Lemma 3.3} Under the same assumption as in Lemma 3.1,
then $b=0$ or $b=1$.
\endproclaim
\demo{\smc Proof} Keep the notations above.
If $b>0$, then $\phi_2$ is actually a morphism.
Thus we have the following commutative diagram:
$$\CD
X   @>h_2>>    C\\
@|            @VVs_2V\\
X   @>>\phi_2>   W_2
\endCD$$
Let $\Ce{0}$ be a saturated subbundle of
 $f_*(\omega_X^{\otimes 2})$ which is generated by $H^0(C,f_*(\omega_X^{\otimes 2}))$.
Because $|2\KX|$ is composed of a pencil and
$\phi_2$ factors through $h_2$, 
$\Ce{0}$ must be a subbundle of rank 1. Let
 ${\Cal E}=f_*(\omega_X^{\otimes 2})$,
we have the following exact sequence
$$0\rw\Ce{0}\rw{\Cal E}\rw\Ce{1}\rw 0$$
and
$$f_*(\omega_{X/C}^{\otimes 2})\rw\Ce{1}\otimes\omega_C^{\otimes -2}\rw 0.$$
Let $ r=rk{\Cal E}=h^0(2K_F)=K_{F_0}^2+\chi(\Co{F_0})\ge 2$. 
By Kawamata's result (\cite{11}), $f_*(\omega_{X/C}^{\otimes 2})$ is semi-positive.
Therefore $\Ce{1}\otimes\omega_C^{\otimes -2}$, as a quotient,
satisfies $\deg(\Ce{1}\otimes\omega_C^{\otimes -2})\ge 0$, i.e.,
$\deg\Ce{1}\ge 4(r-1)(b-1)$. We have
$$\align
h^1(\Ce{0})&\ge h^0(\Ce{1})\ge \deg\Ce{1}+(r-1)(1-b)\\
&\ge 3(r-1)(b-1).
\endalign$$
Noting that $deg\Ce{0}>0$, if $h^1(\Ce{0})>0$, then by  Clifford's theorem,
$$\deg\Ce{0}\ge 2h^0(\Ce{0})-2>h^0(\Ce{0})$$
where $h^0(\Ce{0})=p(2)(X)\ge 4$. We have
$$h^1(\Ce{0})=(h^0(\Ce{0})-deg\Ce{0})+(b-1)<b-1$$
thus $3(r-1)(b-1)<b-1$, which is impossible.
Therefore $h^1(\Ce{0})=0$ and $b=1$.
\qed\enddemo

\proclaim{\smc Lemma 3.4} Under the same assumption as in Lemma 3.1,
we have $p_g(F)\ge 1$.
\endproclaim
\demo{\smc Proof} If $p_g(F)=0$, because $F$ is a surface of general type,
 $q(F)=0$. Therefore  $R^1{h_2}_*\omega_{X'}=0$. By basic formula,
we have
$q(X)=q(X')=b$ and $h^2(\Co{X})=h^2(\Co{X'})=0$. 
If $p_g(X)\ge 1$, we know that $p_g(F)\ge 1$, therefore, under the above assumption,
we must have
$p_g(X)=0$. From Lemma 3.3,
$$\chi(\Co{X})=1-q(X)=1-b\ge 0.$$
which is impossible, because
 $\chi(\Co{X})<0.$
\qed\enddemo

\proclaim{\smc Theorem 3.1} Let $X$ be a nonsingular projective minimal threefold
of general type, suppose that 
$|2\KX|$ be composed of a pencil of surfaces, then 
 $X$ must be of one of the following types:

\noindent (1) $q(F)=0$, $ 1\le K_{F_0}^2\le 3$:
\roster
\item"(11)" $b=1$, $p_g(F)=q(X)=1$, $h^2(\Co{X})=0$, $p_g(X)\ge 2$;
\item"(12)" $b=1$, $1\le p_g(F)\le 3$, $q(X)=1$, $h^2(\Co{X})=0$, $p_g(X)=1$,
   $\chi(\Co{X})=-1$;
\item"(13)" $b=0$, $p_g(F)=1$, $q(X)=h^2(\Co{X})=0$, $p_g(X)\ge 2$. 
\endroster

\noindent (2) $p_g(F)=q(F)=1$, $K_{F_0}^2=2,3$:
\roster
\item"(21)" $b=1$, $q(X)=2$, $h^2(\Co{X})=1$, $p_g(X)\ge 1$;
\item"(22)" $b=1$, $q(X)=1$, $h^2(\Co{X})=0$, $p_g(X)=1$;
\item"(23)" $b=1$, $q(X)=1$, $p_g(X)\ge 2$;
\item"(24)" $b=0$, $q(X)=1$, $h^2(\Co{X})=0$, $p_g(X)\ge 1$;
\item"(25)" $b=0$, $q(X)=0$, $p_g(X)\ge 2$.
\endroster
\endproclaim
\demo{\smc Proof}
From Lemma 3.2, we know that 
 $F$ is of two cases:
(1) $q(F)=0$; (2) $p_g(F)=q(F)=1$.
\vskip0.3cm

\noindent {\smc Case} (1):
 
We have $\triangle_2=degR^1{h_2}_*\omega_{X'/C}
=0$, therefore $q(X)=b$ and 
$$h^2(\Co{X})=h^1({h_2}_*\omega_{X'}).$$

Case$(1)^1$:
 $p_g(X)\ge 2$. It is obvious that 
$|\KX|$ is composed of a pencil of surfaces and
 $\phi_1$ generically factors through $\phi_2$. Take a common birational
modification $f:X'\rw X$ such that  $g_i=\phi_i\circ f$ $(i=1,2)$
is a morphism. We have the following commutative diagram:
$$\CD
X'   @>h>>    C\\
@|            @VVs_2V\\
X'   @>>g_2>   W_2\\
@|             @VVs_1V\\
X'  @>>g_1>    W_1\\
@VfVV            @.\\
X    @.
\endCD$$
Let $g_2:=s_2\circ h$ is a Stein-factorization of
 $g_2$, then $g_1=(s_1\circ s_2)\circ h$ is a Stein-factorization of
 $g_1$. Let $H_1, H_2$ be the general hyperplane section of $W_1, W_2$, 
respectively. We have
 $g_1^*(H_1)\simlin \sum_{i=1}^{a_1}F_i$, $F_i$ is a  fiber of $h$ for every $i$
and $a_1\ge p_g(X')-1$.

If $b=1$, then  $\phi_1$, $\phi_2$ are morphisms.
We may suppose that $X=X'$.
We also have $q(X)=1$. Using a similar method to that in the proof of Lemma 3.3,
one has $h^1(h_*\omega_X)=0$, therefore
$h^2(\Co{X})=0$. Upon an open Zariski subset of $C$,
we have the following exact sequence:
$$0\rw\Co{X'}(K_{X'})\rw\Co{X'}(K_{X'}+g_1^*(H_1))\rw
\oplus_{i=1}^{a_1}\Co{F_i}(K_{F_i})\rw 0.\tag3.1$$
We have the surjective map
$$H^0(K_{X'}+g_1^*(H_1))\rw \oplus_{i=1}^{a_1}H^0(K_{F_i}),$$
thus $p_g(F)=1$, otherwise because
$$\Fi{K_{X'}+g_1^*(H_1)}|_{F_i}=\Fi{K_{F_i}},$$
$\text{dim}\phi_2(X)\ge 2$, a contraction to
our assumption.
Thus  $X$ corresponds to type (11) of the Theorem.

If $b=0$, then $q(X)=0$. Because $\chi(\Co{X})=1+h^2(\Co{X})-p_g(X)<0$, 
$$h^2(\Co{X})\le p_g(X)-2.\tag3.2$$ 
Noting that $|K_{X'}+g_1^*(H_1)|$ is also composed of a pencil of surfaces,
we can easily see that
 $h^0(K_{X'}+g_1^*(H_1))=2p_g(X')-1=2p_g(X)-1$.
 $g_1^*(H_1)\simnum a_1F$, where $a_1=p_g(X')-1$.
From (3.1) and (3.2), we obtain
$$a_1p_g(F)\le p_g(X)-1+h^2(\Co{X})\le 2p_g(X)-3$$
i.e.,  $(p_g(X)-1)p_g(F)\le 2p_g(X)-3.$ Therefore $p_g(F)=1$, 
and then $h_*\omega_{X'}$ is a rank one vector bundle.
 Because $\deg h_*\omega_{X'}>0$,  $h^2(\Co{X})=h^1(h_*\omega_{X'})=0$.
Therefore $X$ corresponds to type (13).

Case$(1)^2$: 
 $p_g(X)\le 1$. From $\chi(\Co{X})=1-q(X)+h^2(\Co{X})-p_g(X)<0$, we get
$q(X)>0$ and then  $b=q(X)=1$, $h^2(\Co{X})=0$, $p_g(X)=1$,  $\chi(\Co{X})=-1$.
$X$ corresponds to type (12).

\vskip0.3cm

\noindent {\smc Case} (2):

In this case, $R^1h_*\omega_{X'}$ is a rank one vector bundle.
Because $R^1h_*\omega_{X'/C}$ is semi-positive,
  $h^1(R^1h_*\omega_{X'})\le 1$.
Note that $h_*\omega_{X'}$ is also a rank one vector bundle and
$b=0,1$. From Riemann-Roch,
we have  $h^1(h_*\omega_{X'})=0$
if $p_g(X)\ge 2$.

Case$(2)^1$:
 $p_g(X)\ge 2$. If $h^1(R^1h_*\omega_{X'})=1$, then $R^1h_*\omega_{X'}
\cong \omega_C$. When $b=1$, then
 $q(X)=2$, $h^2(\Co{X})=1$. $X$ corresponds to type (21);
when $b=0$, then $q(X)=1$, $h^2(\Co{X})=0$. $X$ corresponds to type (24).

If $h^1(R^1h_*\omega_{X'})=0$, then $q(X)=b$. When $b=1$, $X$ corresponds
to type (23);
when $b=0$, $X$ corresponds to type (25).

Case$(2)^2$:
 $p_g(X)\le 1$. From $\chi(\Co{X})<0$, we get $q(X)>0$. $q(X)=b+h^1(R^1h_*\omega_{X'})$.
When $b=0$, then $h^1(R^1h_*\omega_{X'})=1$, $R^1h_*\omega_{X'}\cong\omega_C$.
In this case, $q(X)=p_g(X)=1$ and $h^2(\Co{X})=0$, $\chi(\Co{X})=-1$, $X$ 
corresponds to type (24). When $b=1$, then there is only two possibilities,
i.e., $(q(X),h^2(\Co{X}),p_g(X))=(2,1,1)$
and $(1,0,1)$. The former corresponds to type (21), the latter to type (22).
\qed\enddemo

\proclaim{\smc Corollary 3.1} Let $X$ be a nonsingular minimal projective
threefold of general type,
 if $|2K_X|$ is composed of a pencil of surfaces,
then  $q(X)\le 2$ and $p_g(X)\ge 1$.
\endproclaim

\head 4. Proof of theorem 2\endhead
In this section, we mainly discuss the case when $p_g(X)\le 1$ and always suppose
$|2K_X|$ be composed of a pencil of surfaces. From Theorem 3.1, we see that
$X$ corresponds to type (12), type (21), type (22) and type (24). We keep
the same notations and use the first commutative diagram of the former 
section.

\proclaim{\smc Theorem 4.1} Let $X$ be a nonsingular projective threefold
 with nef and big canonical divisor $K_X$. Suppose $|2K_X|$ be composed of 
a pencil of surfaces, $X$ not corresponding to type (12), then $\fei{5}$ 
is a birational map onto its image.
\endproclaim
\demo{\smc Proof}
Considering the system $|K_{X'}+2f_2^*(K_X)+g_2^*(H_2)|$, we can take a 
standard argument to this situation. Simply, we get from Lemma 3.1 that,
for a general fiber $F$ of $h_2$,
$$\Fi{K_{X'}+2f_2^*(K_X)+g_2^*(H_2)}|_F=\Fi{K_F+2\pi^*(K_{F_0})}=
\Fi{3K_F}.$$
The only exception to the birationality of the 5-canonical map for a 
minimal surface $F_0$ is one with
$$(K_{F_0}^2, p_g(F_0))=(1,2)\ \  \text{\rm or}\ \  (2,3).$$ 
Which just corresponds
to type (12).
\qed\enddemo

\proclaim{\smc Theorem 4.2} Let $X$ be a nonsingular projective threefold
with nef and big canonical divisor $K_X$. Suppose $|2K_X|$ be composed of a 
pencil of surfaces and  $X$ corresponding to
type (12), then $\fei{5}$ is also a birational map onto its image.
\endproclaim
\demo{\smc Proof}
Using the first commutative diagram in \S3, we have
$f_2^*(2K_X)\simlin g_2^*(H_2)+Z_2',$ where $Z_2'$ is the fixed part. 
Take some hyperplane section $\overline{H_2}$ such that $g_2^*(\overline{H_2})
=\sum_{i=1}^{a_2}F_i,$ where $a_2=p(2)\ge 4$ noting that $X$ corresponds
to type $(12)$. At first, we can modify $f_2$ such that
$\sum_{i=1}^{a_2}F_i+Z_2'$ has support with only normal crossings.

Let $D\in |f_2^*(K_X)|$ be the unique effective divisor. Because
$2D\simlin 2f_2^*(K_X),$ there is a hyperplane section $H_2^0$ of $W_2$
in ${\Bbb P}^{p(2)-1}$ such that 
$2D=g_2^*(H_2^0)+Z_2'.$ Set $Z_2':=Z_V+2Z_H,$ where $Z_V$ is the 
vertical part with respect to fibration $h_2:X'\rightarrow C$ and $2Z_H$
is the horizontal part. Thus 
$$D=\frac{1}{2} [g_2^*(H_2^0)+Z_V]+Z_H.$$
Noting that $D$ is a divisor, for a general fiber $F$, 
$Z_H|_F=D|_F\simlin \pi^*(K_{F_0})$ by lemma 3.1.

Considering the ${\Bbb Q}-$divisor 
$$K_{X'}+4f_2^*(K_X)-F-\frac{1}{4}(F_5+\cdots+F_{a_2})-\frac{1}{4}Z_V-
\frac{1}{2}Z_H,$$
set 
$$G:=4f_2^*(K_X)-\frac{1}{4}(F_5+\cdots+F_{a_2})-\frac{1}{4}Z_V-
\frac{1}{2}Z_H$$
and
$$D_0:=\lceil G\rceil=3f_2^*(K_X)+\lceil\frac{Z_H}{2}\rceil-\text{vertical divisors}.$$
For a general fiber $F$, $G-F\simnum\frac{7}{2}f_2^*(K_X)$ is nef and big.
Therefore, by vanishing theorem, $H^1(X', K_{X'}+D_0-F)=0$. 
We then have the surjective map
$$H^0(X', K_{X'}+D_0)\longrightarrow H^0(F, K_F+3\pi^*(K_{F_0})+\lceil
\frac{\pi^*(K_{F_0})}{2}\rceil ).$$
If $F$ is not a surface with $(K^2, p_g)=(1,2)$, then 
$\Fi{K_F+3\pi^*(K_{F_0})+\lceil\frac{\pi^*(K_{F_0})}{2}\rceil}$ 
is birational on $F$. Otherwise, we have the same
statement by Lemma 1.1.
Therefore $\Fi{K_{X'}+D_0}$ is birational and so is $\Fi{5K_{X'}}.$
\qed\enddemo

\head References\endhead
\roster
\item"[1]" Barth, W., Peter, C. and Van de Van, A.: Compact Complex
Surface, Berlin Heidelberg New York Tokyo: Springer 1980.
\item"[2]" Benveniste, X.: Sur les applications pluricanoniques des
	vari\'{e}t\'{e}s de type tr\`{e}s g\'{e}n\'{e}ral en
     dimension 3, {\it Amer.  J. Math.} {\bf 108} (1986), 433--449.
\item"[3]" Bombieri, E.: Canonical models of surfaces of general type,
{\it Publ. Math. Inst. Hautes Etude. Sci.} {\bf 42}, 171-219(1973).
\item"[4]" Chen, M.: An extension of Benveniste-Matsuki's method on 
	6-canonical maps for threefolds, {\it Comm. in Algebra} {\bf
22}, 5759--5767(1994).
\item"[5]" Chen, M.: On pluricanonical maps for threefolds of general type,
{\it J. Math. Soc. Japan}, Vol.{\bf 50}, No.4(1998), to appear.
\item"[6]" Chen, M.: Complex varieties of general type whose canonical
systems are composed with pencils, {\it J. Math. Soc. Japan}, 
Vol.{\bf 51}, No.2(1999), to appear.
\item"[7]" Debarre, O.: Addendum, in\'egaliti\'es num\'eriques pour les surfaces
de type g\'en\'eral, {\it Bull. Soc. Math. France} {\bf 111}/4, 301-302(1983).
\item"[8]" Hironaka, H.: Resolution of singularities of an algebraic
	variety over a field of characteristic zero, {\it Ann. of
Math.} {\bf 79}, 109--326(1964).
\item"[9]" Kawamata, Y.: Cone of curves of algebraic
varieties, {\it Ann. of Math.} {\bf 119}, 603--633(1984).
\item"[10]" Kawamata, Y.: A generalization of Kodaira-Ramanujam's
vanishing theorem, {\it Math. Ann.}, {\bf 261}(1982), 43-46.
\item"[11]" Kawamata, Y.: Kodaira dimension of algebraic fiber spaces over
curves, {\it Invent. Math.} {\bf 66}, 57-71(1982).
\item"[12]" Koll\'{a}r, J.: Higher direct images of dualizing
sheaves I, {\it Ann. of Math.} {\bf 123}, 11--42(1986).
\item"[13]"  Koll\'ar, J.: Higher direct images of dualizing
sheaves II, {\it Ann. of Math.} {\bf 124}, 171--202(1986).
\item"[14]" Lee, S.: Remarks on pluricanonical and the adjoint linear
series, preprint.
\item"[15]" Matsuki, K.: On pluricanonical maps for 3-folds of general
	type, {\it J. Math. Soc. Japan} {\bf 38}, 339--359(1986).
\item"[16]" Miyaoka, Y.: The Chern classes and Kodaira dimension of 
	a minimal variety,
In: {\it Algebraic Geometry}, Sendai, 1985 (Avid. Stud. in Pure Math.
{\bf 10}, 1987, pp. 449--476).
\item"[17]" Ohno, K.: Some inequalities for minimal fibrations of surfaces
of general type over curves, {\it J. Math. Soc. Japan,} {\bf 44}(1992), 643-666.
\item"[18]" Reider, I.: Vector bundles of rank 2 and linear systems
on algebraic surfaces, {\it Ann. of Math.}, {\bf 127}(1988), 309-316.
\item"[19]" Tankeev, S. G.: On n-dimensional canonically polarized varieties and 
varieties of fundamental type, {\it Izv. Akad}, Nauk. SSSR, Ser. Mat.,
{\bf 35}(1971), 31-44; {\it Math. USSR Izv.}, {\bf 5}(1971), 29-43.
\item"[20]" Viehweg, E.: Vanishing theorems, {\it J. reine angew. Math.}, 
{\bf 335}(1982), 1-8.
\endroster

\enddocument